\documentclass[twocolumn]{svjour3u} 
\usepackage[english]{babel}
\usepackage[T1]{fontenc}
\usepackage[utf8]{inputenc}
\usepackage{mathptmx}

\usepackage{amsmath,amsfonts,amsthm,amssymb}
\usepackage{tikz,pgfplots}
\usepgfplotslibrary{groupplots}
\usepackage{tkz-fct}
\usetkzobj{all}
\usetikzlibrary{calc}
\usetikzlibrary{positioning}

\usepackage{epstopdf}
\usepackage{graphicx}
\usepackage{gensymb}
\usepackage{tabularx,booktabs}
\usepackage[numbers]{natbib}

\usepackage{bm}
\newcommand{\bs}{\textbf}
\newcommand{\is}{\textit}

\journalname{Computational Mechanics}

\begin{document}

\title{Orthotropic rotation-free thin shell elements}
\author{Gautam Munglani \and Roman Vetter \and Falk K. Wittel \and Hans J. Herrmann}
\institute{Gautam Munglani \at Computational Physics for Engineering Materials,\\
Institute for Building Materials, ETH Z\"urich,\\
Stefano-Franscini-Platz 3, CH-8093 Z\"urich, Switzerland.\\
\email{gmunglani@ethz.ch}}
\maketitle

\begin{abstract}
A method to simulate orthotropic behaviour in thin shell finite elements is proposed. The approach is based on the transformation of shape function derivatives, resulting in a new orthogonal basis aligned to a specified preferred direction for all elements. This transformation is carried out solely in the undeformed state leaving minimal additional impact on the computational effort expended to simulate orthotropic materials compared to isotropic, resulting in a straightforward and highly efficient implementation. This method is implemented for rotation-free triangular shells using the finite element framework built on the Kirchhoff--Love theory employing subdivision surfaces. The accuracy of this approach is demonstrated using the deformation of a pinched hemispherical shell (with a 18\degree \, hole) standard benchmark. To showcase the efficiency of this implementation, the wrinkling of orthotropic sheets under shear displacement is analyzed. It is found that orthotropic subdivision shells are able to capture the wrinkling behavior of sheets accurately for coarse meshes without the use of an additional wrinkling model. 
\end{abstract}

\keywords{Finite elements \and Rotation-free shells \and Orthotropic materials \and Subdivision surfaces \and Wrinkling} 

\section{Introduction} \label{intro}
Shells possess a unique curved shape which allows them to carry transversal loads primarily by in-plane forces \cite{Ramm2004}. This feature permits them to be extremely slender if necessary, resulting in a wide range of potential applications. Efforts to understand their mechanical capabilities have therefore been crucial in the development of load carrying structures. In conjunction with the finite element method (FEM), shells can be numerically simulated to investigate how mem\-brane-like structures manage different loading conditions, thereby greatly assisting in their design \cite{Bischoff2004}.
 
Biological structures like plant cell walls and skin often possess very thin sheets which undergo large deformations \cite{Farsad2003} that are difficult to capture using traditional simulation methods. These sheets are generally heterogenous in nature with anisotropic material properties \cite{Gibson2012}. Some even have preferential orientations due to the presence of crosslinked polymers, giving them orthotropic or transversely isotropic material properties \cite{Burgert2006}. Due to their geometry, these sheets are often modelled with computationally efficient thin shell finite elements \cite{Holzapfel1996}. It is therefore essential that the type of shell elements used to approximate these structures are capable of simulating non-isotropic material properties. 

A traditional way to treat thin shells numerically is with the use of the Kirchhoff--Love theory. This theory uses a set of assumptions to effectively capture the properties of geometrically exact thin shells in terms of its middle surface. These assumptions are that the shell is sufficiently thin, and that straight material lines initially normal to the middle surface of the shell remain unstretched, straight and normal after deformation such that the normal stress in the thickness direction is negligible. Shells simulated with this approach have been shown to be quite accurate on a variety of numerical benchmarks \cite{Eugenio2000,Bischoff2004}. 

The choice of shape functions is crucial in the analysis of thin shells of Kirchhoff--Love type due to the well-researched $C^1$ continuity requirement \cite{Zienkiewicz2000}. This requirement is due to the presence of second-order derivatives of displacement, which leads to a fourth-order equilibrium equation. This in turn calls for continuous first order derivatives across element boundaries. In order to tackle the continuity requirement, higher-order conforming shape functions with additional degrees of freedom in Hsieh-Clough-Tocher triangles \cite{Clough1965} and Hermite quadrilaterals \cite{Zienkiewicz2000} have previously been used in shell elements. Alternatively, this requirement has been ignored and non-conforming elements with $C^0$ continuity have been combined with other assumptions \cite{Reddy2004}. Recently, isogeometric analysis using NURBS \cite{Hughes2005}, as well as subdivision surfaces \cite{Cirak2000,Cirak2001} have also been developed to satisfy these requirements. Subdivision surfaces in particular have been implemented in triangular shells with a rotation-free formulation, i.e. requiring only displacement degrees of freedom at the mesh nodes. This rotation-free formulation greatly saves on computational cost due to the significant reduction in the number of degrees of freedom, thereby gaining prominence in recent years \cite{Eugenio2000,Linhard2007,Onate2005,Brunet2006}. Due to their formulation, triangular rotation-free elements like subdivision shells pose a unique set of challenges to the simulation of orthotropic material properties required to model deforming sheets \cite{Valdes2007,Valdes2009}. 
 
In this article, we present a straightforward method to simulate orthotropy in Kirchhoff--Love rotation-free shell elements. This method is based on the transformation of the directional derivatives of the element shape functions. It is computationally efficient, independent of the choice of shape functions, and can be coupled to a wide range of problems without further alterations. We use Loop subdivision shell elements \cite{Loop1987} to demonstrate its accuracy for geometrically nonlinear problems, which include the pinched hemisphere benchmark and the problem of thin sheets wrinkling under applied shear. 

This article is divided into the following sections: first, the typical continuum formulation for geometrically exact Kirchhoff--Love thin shells is described. Section 3 then provides a short overview on the development of the basis orientation followed by a detailed explanation of the transformation method. This is followed by a few words on the constitutive law and discretization. Finally, the accuracy of the method using the pinched hemisphere benchmark is demonstrated which culminates with the application of orthotropic subdivision shells to the analysis of thin sheets subjected to shear. 

\section{Kinematics} \label{sec:kine}
To begin, the standard formulation of geometrically exact Kirchhoff--Love thin shells is revisited. The deformation of the thin shell is constructed using the classical stress-resultant formulation in curvilinear coordinates~\cite{Simo1989}. This formulation describes the shell in terms of its middle surface by integrating the stress through its thickness. 

We begin by introducing some standard notation. The reference geometry of the shell is characterized by a middle surface $\overline{\Omega} \subset \mathbb{E}^3$, with boundary $\overline{\Gamma} = \partial\overline{\Omega}$ and thickness \is{h} which is considered to be small compared to the planar dimensions of the shell. With the action of applied loads, the shell deforms into a new configuration characterized by an altered middle surface $\Omega \subset \mathbb{E}^3$. The position vectors of an arbitrary material point in the shell, denoted by $\mathbf{\overline{r}}$ and $\mathbf{r}$ in the reference and deformed configurations respectively, are parametrized in terms of their curvilinear coordinates $\{\theta^1, \theta^2, \theta^3\}$ as
\begin{equation} \label{eq:rx}
\begin{split}
\overline{\bs{r}} (\theta^1, \theta^2, \theta^3) &= \overline{\bs{x}} (\theta^1, \theta^2) + \theta^3 \, \overline{\bs{a}}_{3} \, (\theta^1, \theta^2) \, , \\
\bs{r} (\theta^1, \theta^2, \theta^3) &= \bs{x} (\theta^1, \theta^2) + \theta^3 \, \bs{a}_{3} \, (\theta^1, \theta^2) \, ,
\end{split}
\end{equation}
where $\theta^3 \in [-\is{h}/2, \, \is{h}/2]$ represents the thickness coordinate, $\overline{\bs{x}}$ and $\bs{x}$ are the parametric representations of the middle surfaces of the configurations $\overline{\Omega}$ and $\Omega$ respectively, while $\mathbf{\overline{a}}_{3}$ and $\mathbf{a}_{3}$ are the shell director vectors in the reference and deformed configurations respectively. The change in the material point from $\overline{\bs{x}}$ to $\bs{x}$ is described in terms of a displacement field $\bs{u} = \bs{x} - \overline{\bs{x}}$. This parametrization allows the description of the middle surface coordinate basis which span the tangent space $T$ as 
\begin{equation} \label{eq:basis}
\overline{\bs{a}}_\alpha = \dfrac{\partial \overline{\bs{x}}}{\,\, \partial \theta^\alpha} = \overline{\bs{x}}_{,\, \alpha} \,\, , 
\hspace{15 pt}
\bs{a}_\alpha = \dfrac{\partial \bs{x}}{\,\, \partial \theta^\alpha} = \bs{x}_{,\,\alpha} \, ,
\end{equation}
where Greek indices indicate integers 1 and 2, and a comma indicates partial differentiation. In addition, it should be noted that henceforth, Einstein summation is used and Latin indices indicate integers from 1 to 3. By applying the Kirchhoff constraints, the shell directors can be described as
\begin{equation} \label{eq:director}
\begin{split}
\overline{\bs{a}}_3 &= \frac{\overline{\bs{a}}_1 \, \times \, \overline{\bs{a}}_2}{\overline{J}} \, ,
\hspace{15 pt}
\overline{J} = \| \overline{\bs{a}}_{1} \times \overline{\bs{a}}_{2} \| \, \\
\bs{a}_3 &= \frac{\bs{a}_1 \, \times \, \bs{a}_2}{J} \, ,
\hspace{15 pt}
J = \| \bs{a}_{1} \times \bs{a}_{2} \| \, \\
\end{split}
\end{equation}
where \is{J} is the Jacobian determinant. The infinitesimal area element is now expressed as $\mathrm{d}\overline{\Omega} = \overline{J} \, \mathrm{d}\theta^1 \mathrm{d}\theta^2$. By differentiating Eq.(\ref{eq:rx}), the covariant basis of an arbitrary point in the shell can be described as
\begin{equation} \label{eq:diff}
\begin{split}
\frac{\partial \overline{\bs{r}}}{\partial \theta^\alpha} &= \frac{\partial \overline{\bs{x}}}{\partial \theta^\alpha} + \theta^3 \, \frac{\overline{\bs{a}}_3} {\partial \theta^\alpha}  \, ,
\hspace{25 pt}
\frac{\partial \overline{\bs{r}}}{\partial \theta^3} = \frac{\partial \overline{\bs{x}}}{\partial \theta^3} \, ,\\
\frac{\partial \bs{r}}{\partial \theta^\alpha} &= \frac{\partial \bs{x}}{\partial \theta^\alpha} + \theta^3 \, \frac{\bs{a}_3} {\partial \theta^\alpha}  \, , 
\hspace{25 pt}
\frac{\partial \bs{r}}{\partial \theta^3} = \frac{\partial \bs{x}}{\partial \theta^3} \, .\\
\end{split}
\end{equation} 
Substituting Eq.(\ref{eq:basis}) in the above yields
\begin{equation} \label{eq:gbasis}
\begin{split}
\overline{\bs{g}}_{\alpha} &= \overline{\bs{a}}_{\alpha} + \theta^3 \, \overline{\bs{a}}_{3, \,\alpha} \, , 
\hspace{30 pt}
\overline{\bs{g}}_3 = \overline{\bs{a}}_3 \, ,\\
\bs{g}_{\alpha} &= \bs{a}_{\alpha} + \theta^3 \, \bs{a}_{3, \, \alpha} \, , 
\hspace{30 pt}
\bs{g}_3 = \bs{a}_3 \, ,\\
\end{split}
\end{equation}  
while the corresponding covariant coefficients of the metric tensors are
\begin{equation} \label{eq:gmetric}
\overline{g}_{i j} = \overline{\bs{g}}_{i}  \, \cdot \, \overline{\bs{g}}_{j} \, ,
\hspace{15 pt}
g_{i j} = \bs{g}_{i}  \, \cdot \, \bs{g}_{j} \, .
\end{equation}
The covariant coefficients of the surface metric tensor are 
\begin{equation} \label{eq:ametric}
\overline{a}_{\alpha \beta} = \overline{\bs{a}}_{\alpha}  \, \cdot \, \overline{\bs{a}}_{\beta} \, ,
\hspace{15 pt}
a_{\alpha \beta} = \bs{a}_{\alpha}  \, \cdot \, \bs{a}_{\beta} \, ,
\end{equation}
while the covariant coefficients of the shape tensor are
\begin{equation} \label{eq:sff}
\begin{split}
\overline{b}_{\alpha \beta} &= - \overline{\bs{a}}_{3, \,\alpha} \, \cdot \, \overline{\bs{a}}_{\beta} = \overline{\bs{a}}_{3} \, \cdot \, \overline{\bs{a}}_{\alpha, \, \beta} \, , \\
b_{\alpha \beta} &= - \bs{a}_{3, \,\alpha} \, \cdot \, \bs{a}_{\beta} = \bs{a}_{3} \, \cdot \, \bs{a}_{\alpha, \, \beta} \, .\\
\end{split}
\end{equation} 
This formulation uses the Green-Lagrange strain ($E_{ij}$), whose coefficients in curvilinear coordinates are defined as
\begin{equation} \label{eq:glstrain1}
E_{ij} = \frac{1}{2} (g_{ij} - \overline{g}_{ij}) \, .
\end{equation}
Using additive decomposition, the Green-Lagrange strain tensor can be split into the non-zero coefficients of the membrane $\alpha_{\alpha \beta}$ and bending $\beta_{\alpha \beta}$ strain tensors. Neglecting higher order terms ($\mathcal{O}(\theta^3)^2$) and considering Eq.(\ref{eq:gbasis}) and (\ref{eq:gmetric}), we arrive at   
\begin{equation} \label{eq:glstrain3}
E_{\alpha \beta} = \alpha_{\alpha \beta} + \theta^3 \beta_{\alpha \beta} \, ,
\end{equation}
where the coefficients of the strain tensors are
\begin{equation} \label{eq:strtensor}
\alpha_{\alpha \beta} = \frac{1}{2}  (a_{\alpha \beta} - \overline{a}_{\alpha \beta}) \, ,\\
\hspace{15 pt}
\beta_{\alpha \beta} = \overline{b}_{\alpha \beta} - b_{\alpha \beta} \, .\\
\end{equation}
This neglection of the higher order terms makes this a first order shell theory, which is valid for small thickness \is{h}.
\section{Shape function derivative transformation} \label{sec:ortho}
\begin{figure}[t] 
\centering
\includegraphics[width=8.4cm]{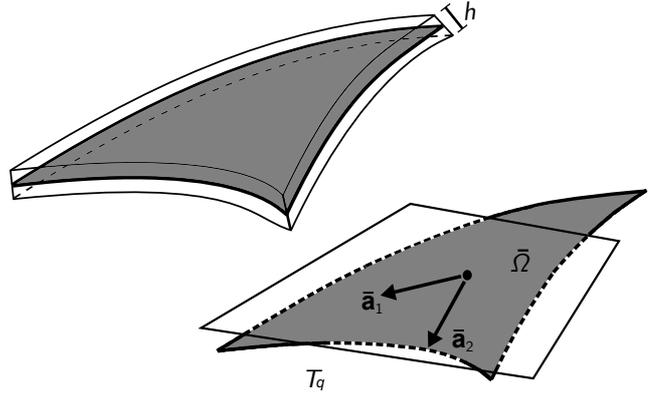}
\caption{Thin shell with reference thickness \is{h} represented by its middle surface $\overline{\Omega}$ with non-orthogonal basis $\overline{\bs{a}}_{\alpha}$ spanning tangent space $T_q$ at point \is{q}} \label{fig:middle}
\end{figure}

Now that the basic notation has been reviewed, we describe how the element basis is traditionally oriented based on the shape functions, followed by our approach to simulate orthotropy. As stated in the previous section, the tangent space $T$ of the middle surface of a shell is spanned by the in-plane basis vectors of the reference and deformed configurations $\overline{\bs{a}}_{\alpha}$ and $\bs{a}_{\alpha}$ respectively. When the continuum is discretized into elements, each quadrature point \is{q} of an element acquires its own tangent space which we shall refer to as $T_q$ as can be seen in Fig.(\ref{fig:middle}). It should be noted that the number of quadrature points per element has no bearing on our proposed method as the transformation occurs on all available $T_q$. Due to this property, the \is{reduced integration} approach of a single quadrature point per element is used henceforth. 

\begin{figure}[h] 
\centering
\includegraphics[width=8.4cm]{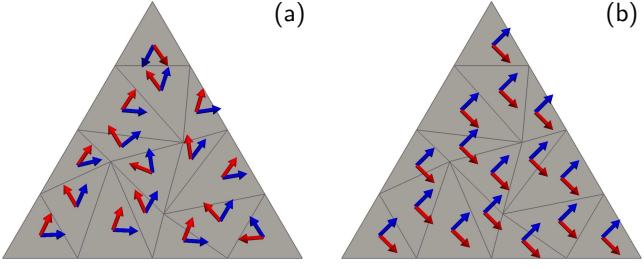} 
\caption{(a) The orientation of the basis is originally non-orthogonal and non-aligned. (b) The transformed basis is orthonormal and aligned to the preferred direction} \label{fig:align}
\end{figure}

\begin{figure}[h]
\centering
\includegraphics[width=8.4cm]{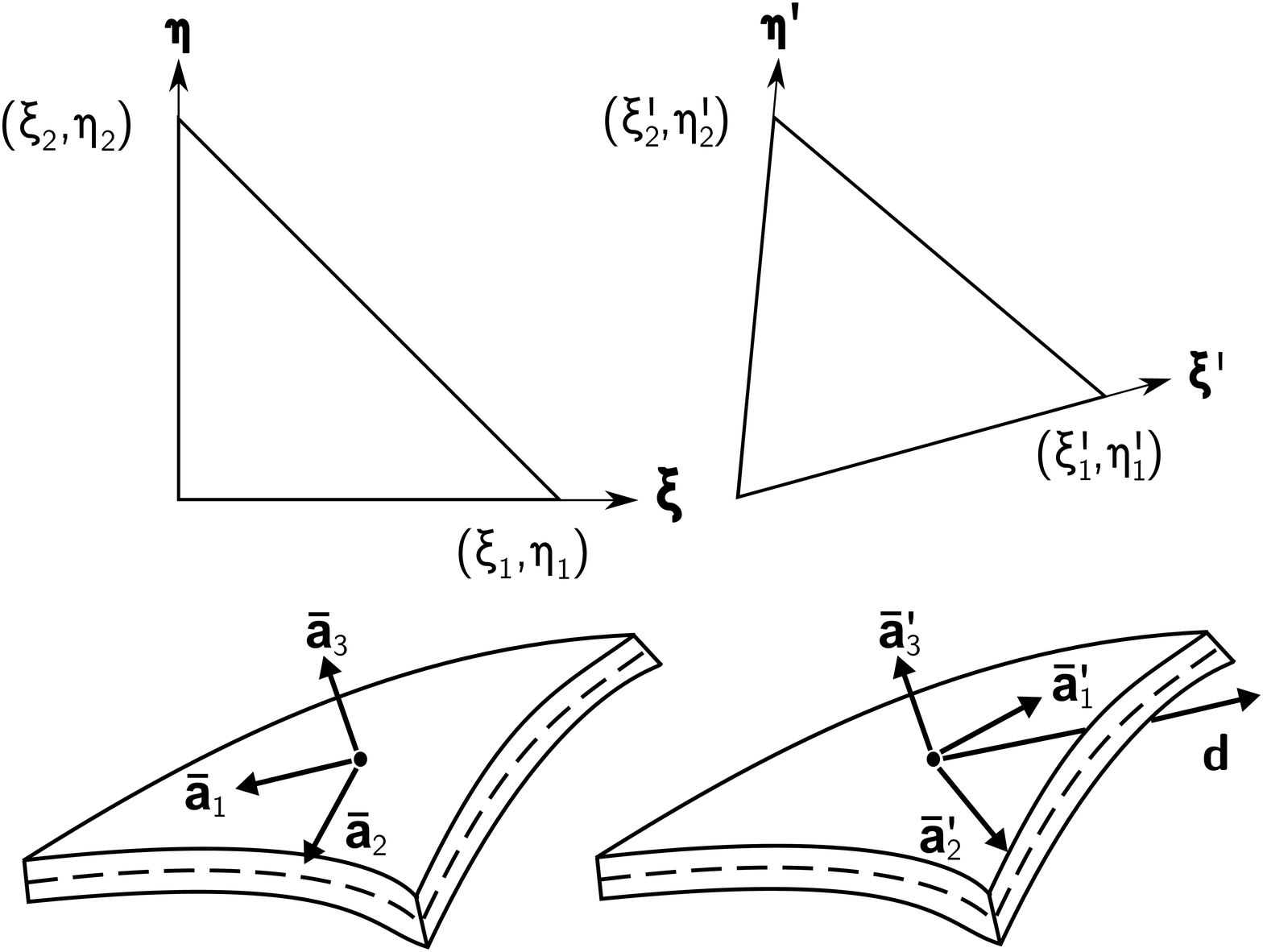} 
\caption{The basis transformation corresponds to a change in the shape of the master triangle from $\bm{\xi}$ and $\bm{\eta}$ to the directions stipulated by $\bm{\xi}'$ and $\bm{\eta}'$} 
\label{fig:stand}
\end{figure}

By construction, triangular elements generally have non-orthogonal basis vectors in $T_q$. This also applies to curved regular meshes, where the orientation of each triangular element varies from its neighbours, resulting in basis vectors that are not aligned between elements (Fig.(\ref{fig:align}a)). Isotropic materials do not depend on the orientation or orthonormality of the basis as their properties remain constant in all directions, unlike the case for orthotropic materials. Therefore, to confer orthotropic properties to a surface, the element basis vectors are transformed to an orthogonal basis and aligned to a preferred direction across the surface as seen in Fig.(\ref{fig:align}b). The discretized form of the basis vectors can be written as a linear combination of shape function derivatives
\begin{equation} \label{eq:qp}
\overline{\bs{a}}_{1} = \sum_{I = 1}^{n_{\mathrm{sf}}} N_{I, \bm{\xi}} \, \overline{\bs{x}}_{I} \, , 
\hspace{15 pt} 
\overline{\bs{a}}_{2} = \sum_{I = 1}^{n_{\mathrm{sf}}} N_{I, \bm{\eta}} \, \overline{\bs{x}}_{I} \, ,
\end{equation}
where $\{N_I, I=1,...,n_{\mathrm{sf}}\}$ are the shape functions and $\overline{\bs{x}}_{I}$ are the nodal positions, with $n_{\mathrm{sf}}$ being the number of mesh nodes accounted for by each quadrature point. These derivatives are defined by directions $\bm{\xi}, \bm{\eta} \in \mathbb{E}^3$ along the sides of the standard triangular master element as shown in Fig.(\ref{fig:stand}). Basis vectors are in fact manifold versions of directional derivatives taken along the sides of this master element. Technically, the basis can be transformed directly using an orthogonalization method like the Gram--Schmidt process \cite{Greub1975}. This would result in an orthonormal basis that spans the same space $T_q$ as the original shown in Fig.(\ref{fig:middle}). However, this basis would be mismatched with the shape function derivatives shown in Eq.(\ref{eq:qp}) and its second derivatives. Alternatively, one could directly transform the covariant strain tensors $\alpha_{\alpha \beta}$ and $\beta_{\alpha \beta}$ in Eq.(\ref{eq:strtensor}) and the contravariant resultant stresses shown later in Sec.(\ref{sec:const}) using the principal directions of orthotropy \cite{Dujc2012}. This procedure would need to be performed at every time step, which can be very computationally costly. Furthermore, the shape function derivatives would now be mismatched with the strains and resultant stresses.  

A straightforward procedure to create a basis matching with the first and second derivatives of its shape functions is to use a transformation that directly alters the first and second derivatives of the shape functions $N_I$ in $\overline{\Omega}$ leading to a new orthogonal basis without further adjustments. By ensuring that the new basis vectors in each $T_q$ are individually aligned to the preferred direction, this method is able to simultaneously transform and align the basis of triangular elements. An advantage of this procedure is that it only needs to be performed once for every element in the mesh for the reference configuration in a Total Lagrangian formulation. Other works \cite{Valdes2007,Valdes2009} have used alternate transformation matrices built using the deformed configuration. This implies that, as with the direct transformation of the strains and resultant stresses, the process of orthogonalization is repeated for every deformed configuration.

The first step in this procedure is to introduce a vector to indicate the preferred direction of the material. This shall be denoted by the unit vector $\bs{d} \in \mathbb{E}^3$, which is not parallel to $\overline{\bs{a}}_3$ anywhere on the middle surface $\overline{\Omega}$. In order to ensure that the reference basis of each element is aligned, $\bs{d}$ is identical for every quadrature point in the mesh. The curvature of the shell surface predicates that $\bs{d}$ will most likely not lie in $T_q$. Therefore, it needs to be projected onto $T_q$ using the relation
\begin{equation} \label{eq:proj}
\hat{\bs{a}}_1 = \bs{d} - (\bs{d} \cdot \overline{\bs{a}}_3) \, \overline{\bs{a}}_3 \, ,
\end{equation}    
where $\hat{\bs{a}}_1$ is the new unscaled in-plane basis vector that points towards the preferred direction while lying in $T_q$. The other new unscaled basis vector $\hat{\bs{a}}_2$ points in the direction perpendicular to the preferred direction, henceforth referred to as the perpendicular direction
\begin{equation} \label{eq:perp}
\hat{\bs{a}}_2 = \overline{\bs{a}}_3 \times \hat{\bs{a}}_1   \, .
\end{equation}

The unscaled basis $\hat{\bs{a}}_{\alpha}$ now needs to be rescaled by the lengths of the original basis vectors $\overline{\bs{a}}_{\alpha}$ to the new in-plane basis in the reference configuration, $\overline{\bs{a}}'_{\alpha}$. This is required as the basis transformation is non-orthogonal, thereby not preserving vector lengths. The angle between the original basis $\overline{\bs{a}}_{1}$ and the new preferentially aligned $\hat{\bs{a}}_{1}$ is given by
\begin{equation} \label{eq:theta}
\theta = \arctan \left(\frac{\|\hat{\bs{a}}_{1} \times \overline{\bs{a}}_{1}\|}{\hat{\bs{a}}_{1} \,\cdot \, \overline{\bs{a}}_{1}}\right) \, .
\end{equation}    
Since $\overline{J}$ has to remain constant, $\overline{\bs{a}}_3$ can be treated as the reference around which the original basis is rotated. Using the relation $\|\overline{\bs{a}}_{3}\| = \|\overline{\bs{a}}_{1}\| \, \|\overline{\bs{a}}_{2}\| \, \sin\theta$, the new orthogonal basis $\overline{\bs{a}}'_{\alpha}$ can be defined as
\begin{equation} \label{eq:rescale}
\overline{\bs{a}}'_{\alpha} = \hat{\bs{a}}_{\alpha} \, \frac{\|\overline{\bs{a}}_{3}\|}{\|\overline{\bs{a}}_{\beta}\| \, \sin\theta} \,
\hspace{15 pt}
\text{where} \,\, \alpha \neq \beta.
\end{equation}
Since $\overline{\bs{a}}'_{\alpha}$ has been constructed, and $\overline{\bs{a}}'_{3} = \overline{\bs{a}}_{3}$, we can look at altering the shape function derivatives, $N_{I, \bm{\xi}}$ and $N_{I, \bm{\eta}}$ respectively. In their original state,  $\bm{\xi} \, \text{and} \, \bm{\eta}$ correspond to the three-dimensional Cartesian unit vectors $\bs{e}_1$ and $\bs{e}_2$ respectively. 
\begin{equation} 
\begin{split}
\bm{\xi} &= [\xi_1 \hspace{10 pt} \xi_2 \hspace{10 pt} \xi_3] = [1\hspace{10 pt} 0\hspace{10 pt} 0] \, , \\
\bm{\eta} &= [\eta_1 \hspace{10 pt} \eta_2 \hspace{10 pt} \eta_3] = [0\hspace{10 pt} 1\hspace{10 pt} 0] \, . \\
\end{split}
\end{equation}
These directions need to be modified to $\bm{\xi}'$ and $\bm{\eta}'$ respectively. These new directions are found by constructing a transformation matrix \bs{T} between $\overline{\bs{a}}'_{k}$ and $\overline{\bs{a}}_{k}$
\begin{equation} \label{eq:transformer}
\bs{T} = (\overline{\bs{a}}_k \otimes \bs{e}_k)^{-1} \, (\overline{\bs{a}}'_k \otimes \bs{e}_k) = [T_{ij}] = 
\begin{bmatrix} 
[T_{\alpha \beta}] & 0 \\ 
0 & 1 \\
\end{bmatrix}
\, ,
\end{equation}  
where $\bs{e}_k$ denotes the in-plane Cartesian unit vectors and $[T_{\alpha \beta}]$ is a non-orthogonal $2 \times 2$ matrix representing an in-plane basis transformation with coefficients
\begin{equation}
\renewcommand\arraystretch{1.6}
[T_{\alpha \beta}] = 
\begin{bmatrix}
\eta'_1 & \xi'_1\\
\eta'_2 & \xi'_2\\
\end{bmatrix}
\, ,
\end{equation}
which describes the shape of the new non-right angled master triangle as can be seen in Fig.(\ref{fig:stand}). Application of $[T_{\alpha \beta}]^T$ modifies $N_{I, \bm{\xi}}$ and $N_{I, \bm{\eta}}$ to derivatives with directions $\bm{\xi}'$ and $\bm{\eta}'$ respectively. This transformation also applies to the second derivatives of $N_I$ if higher-order shape functions are used. This is shown by
\begin{equation}
\begin{split}
\nabla ' (N_I) &= [T_{\alpha \beta}]^\mathrm{T} \, \nabla (N_I) \, ,\\
H'(N_I) &= [T_{\alpha \beta}]^\mathrm{T} \, H(N_I) \, [T_{\alpha \beta}] \, , \\ 
\end{split}
\end{equation}
where $\nabla ' (N_I)$ and $H'(N_I)$ represent the modified gradient and Hessian matrix of the shape functions respectively. In tensor form, this can be rewritten as
\begin{equation} \label{eq:Nfirst}
\renewcommand\arraystretch{1.6}
\begin{bmatrix}
N_{I, \bm{\eta}'} \\ 
N_{I, \bm{\xi}'} \\ 
\end{bmatrix}
=
\begin{bmatrix}
\eta'_1 & \eta'_2 \\
\xi'_1 & \xi'_2 \\
\end{bmatrix}
\begin{bmatrix}
N_{I, \bm{\eta}} \\ 
N_{I, \bm{\xi}} \\ 
\end{bmatrix}
\, ,
\end{equation}
\begin{equation} \label{eq:Nsecond}
\renewcommand\arraystretch{1.6}
\begin{bmatrix}
N_{I, \bm{\eta}' \bm{\eta}'} \\ 
N_{I, \bm{\xi}' \bm{\xi}'} \\ 
N_{I, \bm{\xi}' \bm{\eta}'} \\
N_{I, \bm{\eta}' \bm{\xi}'} \\ 
\end{bmatrix}
=
\begin{bmatrix}
\eta'_1 \, \eta'_1 & \eta'_2 \, \eta'_2 & \eta'_1 \, \eta'_2 & \eta'_1 \, \eta'_2\\ 
\xi'_1 \, \xi'_1 & \xi'_2 \, \xi'_2 & \xi'_1 \, \xi'_2 & \xi'_1 \, \xi'_2\\ 
\eta'_1 \, \xi'_1 & \eta'_2 \, \xi'_2 & \eta'_1 \, \xi'_2 & \eta'_2 \, \xi'_1\\ 
\eta'_1 \, \xi'_1 & \eta'_2 \, \xi'_2 & \eta'_2 \, \xi'_1 & \eta'_1 \, \xi'_2\\ 
\end{bmatrix}
\begin{bmatrix}
N_{I, \bm{\eta} \bm{\eta}} \\ 
N_{I, \bm{\xi} \bm{\xi}} \\ 
N_{I, \bm{\xi} \bm{\eta}} \\
N_{I, \bm{\eta} \bm{\xi}} \\ 
\end{bmatrix}
\, .
\end{equation}

The first and second derivatives of the shape functions have now been transformed, allowing the new orthogonal basis to be constructed by replacing $N_{I, \bm{\xi}}$ and $N_{I, \bm{\eta}}$ with $N_{I, \bm{\xi}'}$ and $N_{I, \bm{\eta}'}$ in Eq.(\ref{eq:qp}). The altered second derivatives of the shape functions given in Eq.(\ref{eq:Nsecond}) are used to build the first derivatives of the basis $\overline{\bs{a}}_{\alpha, \, \beta}$. The coefficients of the surface metric tensor and shape tensor defined in Eq.({\ref{eq:ametric}}-{\ref{eq:sff}}) respectively are then determined using this altered basis. This finally results in updated membrane $\alpha_{\alpha \beta}$ and bending $\beta_{\alpha \beta}$ strain tensors shown in Eq.(\ref{eq:strtensor}). Once these geometric quantities have been built, the constitutive law in curvilinear coordinates needs to be defined for the orthotropic case. This combined with the assembly of the stress resultants with the new orthogonal and aligned basis is described in the next section.  
\section{Constitutive model and discretization} \label{sec:const} 
The material properties of the shell are determined by the St. Venant--Kirchhoff constitutive law. The strain energy density function for the orthotropic case is altered from the Koiter energy density functional \cite{Koiter1970} to
\begin{equation} \label{eq:koiter}
W = \frac{K^{\alpha \beta}}{2} \left[h \, H^{\alpha \beta \gamma \delta} \, \alpha_{\alpha \beta} \, \alpha_{\gamma \delta} + \frac{h^3}{12} \, H^{\alpha \beta \gamma \delta} \, \beta_{\alpha \beta} \, \beta_{\gamma \delta} \right] ,
\end{equation}
where the coefficients of the elasticity tensor $H^{\alpha \beta \gamma \delta}$ are
\begin{equation} \label{eq:etensor}
H^{\alpha \beta \gamma \delta} = \nu_1 \overline{a}^{\alpha \beta} \, \overline{a}^{\gamma \delta} + \frac{(1 - \nu_1)}{2} (\overline{a}^{\alpha \gamma} \overline{a}^{\beta \delta} + \overline{a}^{\alpha \delta} \overline{a}^{\beta \gamma}) \, ,
\end{equation}
and the coefficients of the stiffness $K^{\alpha \beta}$ are
\begin{equation} \label{eq:stiffness}
\left[K^{11} \hspace{10 pt} K^{22} \hspace{10 pt} K^{12} \right]
= 
\left[\frac{E_1}{1-\nu_1 \nu_2} \hspace{10 pt} \frac{E_2}{1-\nu_1 \nu_2} \hspace{10 pt} \frac{G_{12}}{1-\nu_1} \right]. 
\end{equation}
The above relation includes the elastic modulus $E$, the Poisson ratio $\nu$ and the shear modulus $G$, where the index 1 represents the preferred direction \bs{d} of the material, and 2 represents the perpendicular direction. The moduli and Poisson ratios are related by $E_1 \, \nu_2 = E_2 \, \nu_1$. The derivative of the strain energy density functional with respect to the membrane strains $\alpha_{\alpha \beta}$ and bending strains $\beta_{\alpha \beta}$ gives the resultant membrane stresses $n^{\alpha \beta}$ and bending stresses $m^{\alpha \beta}$ of the element respectively
\begin{equation} \label{eq:stressvoight}
\begin{split}
n^{\alpha \beta} &= \frac{\partial W}{\partial \alpha_{\alpha \beta}} = h \, K^{\alpha \beta} \, H^{\alpha \beta \gamma \delta} \, \alpha_{\gamma \delta} \, , \\
m^{\alpha \beta} &= \frac{\partial W}{\partial \beta_{\alpha \beta}} = \frac{h^3}{12} \, K^{\alpha \beta} \, H^{\alpha \beta \gamma \delta} \, \beta_{\gamma \delta} \, . \\
\end{split}
\end{equation} 
Now that the form of the stress resultants has been elucidated, we can proceed to the  discretized form of the equilibrium equations. This is achieved by approximating the minimization of the total potential energy, which is obtained by summing the contribution of the internal elastic energy ($\phi^{\mathrm{int}}$) with the external energy ($\phi^{\mathrm{ext}}$) according to
\begin{equation} \label{eq:totalenergy}
\phi[\bs{u}] = \phi^{\mathrm{ext}} [\bs{u}] + \phi^{\mathrm{int}} [\bs{u}] \, ,
\end{equation}  
The internal energy is the Koiter strain energy density is integrated over the reference middle surface, and the external energy is the sum of the external load \bs{q} per unit surface area and the traction \bs{N} per unit edge length, given as
\begin{equation} \label{eq:intenergy}
\begin{split}
\phi^{\mathrm{int}} [\bs{u}] &= \int_{\overline{\Omega}} W \, \mathrm{d}\overline{\Omega} \, , \\
\phi^{\mathrm{ext}} [\bs{u}] &= -\int_{\overline{\Omega}} \bs{q} \cdot \bs{u} \, \mathrm{d}\overline{\Omega} - \int_{\overline{\Gamma}} \bs{N} \cdot \bs{u} \, \mathrm{d}\overline{\Gamma} \, , \\
\end{split}
\end{equation}
respectively. To solve the elastic energy minimization problem for the displacement field \bs{u}, the first variation of $\phi$ is taken and augmented with the inertial term containing the mass matrix. This results in
\begin{equation} \label{eq:firstvariation}
0 = \bs{f}_I^{\mathrm{int} } - \bs{f}_I^{\mathrm{ext}} + \sum_{J} M_{IJ} \, \ddot{\bs{u}}_{J} \, ,
\end{equation} 
where
\begin{equation} \label{eq:firstvariationsplit}
\begin{split}
&\bs{f}_I^{\mathrm{int}} = -\int_{\overline{\Omega}}  \left(n^{\alpha \beta} \, \dfrac{\partial \alpha_{\alpha \beta}}{\partial \bs{u}_I} + m^{\alpha \beta} \, \dfrac{\partial \beta_{\alpha \beta}}{\partial \bs{u}_I}\right) \, \mathrm{d} \overline{\Omega} \, , \\
&\bs{f}_I^{\mathrm{ext}} = \int_{\overline{\Omega}} \bs{q} \, N_I \, \mathrm{d} \overline{\Omega} + \int_{\overline{\Gamma}} \bs{N} \, N_I \, \mathrm{d}\overline{\Gamma} \, , \\
& M_{IJ} = \int h \, \rho \, N_I  N_J \, \mathrm{d} \overline{\Omega} \, . \\
\end{split}
\end{equation}
In the above relation, $M_{IJ}$ is the mass matrix for dynamic analysis.  Eq.(\ref{eq:firstvariation}) can now be evaluated element-wise using a quadrature rule. The \is{constant-average acceleration method} is used for time integration. This is obtained by setting $\gamma = 1/2$ and $\beta = 1/4$ in the widely used Newmark family of methods \cite{Newmark1959,Reddy2004}. Further details on this setup can be found in Refs.\cite{Cirak2000,Cirak2001}. 

Methods that directly transform the components of $\bs{f}^{\mathrm{int}}$ will require the additional transformation of $\bs{f}^{\mathrm{ext}}$ to the new local coordinate system in order to implement boundary conditions which involve shape function derivatives. All rotation-free shells require this supplementary step for their implementation of certain natural boundary conditions \cite{Onate2005, Green2004,Brunet2006,Linhard2007,Cirak2011,Kiendl2009,Nguyen2015} including frictional contact \cite{Wriggers2006}, in-plane shear traction based conditions \cite{Long2012} and displacement-dependent pressure loads \cite{Schwei1984}. Some widely used formulations additionally require these shape function derivatives for curvature gradient calculations \cite{Onate2005,Linhard2007} or hourglass stabilization \cite{Brunet2006}. Our proposed method does not require these further transformations. 

The details of the simulation method and its advantages have been presented and are now followed by two numerical studies in the next section to demonstrate the efficiency and accuracy of orthotropic subdivision shells. 
\section{Numerical studies} \label{sec:num}
We present two examples to demonstrate the accuracy of the described transformation method. The first is a standard numerical benchmark of orthotropic behaviour for geometrically nonlinear problems. The second showcases the efficiency of orthotropic subdivison shell elements by analyzing wrinkling behavior of sheets. Both examples have been simulated using a single quadrature point per element (\textit{reduced integration}), located at its barycenter. It has been argued in Refs.{\cite{Cirak2000,Cirak2001} that a single quadrature point in adequate for the simulation of geometrically nonlinear large deformation problems with subdivison shells.
\subsection{Pinched hemispherical shell} \label{subsec:pinched}
\begin{figure}[h]
\centering
\includegraphics[width=8.4cm]{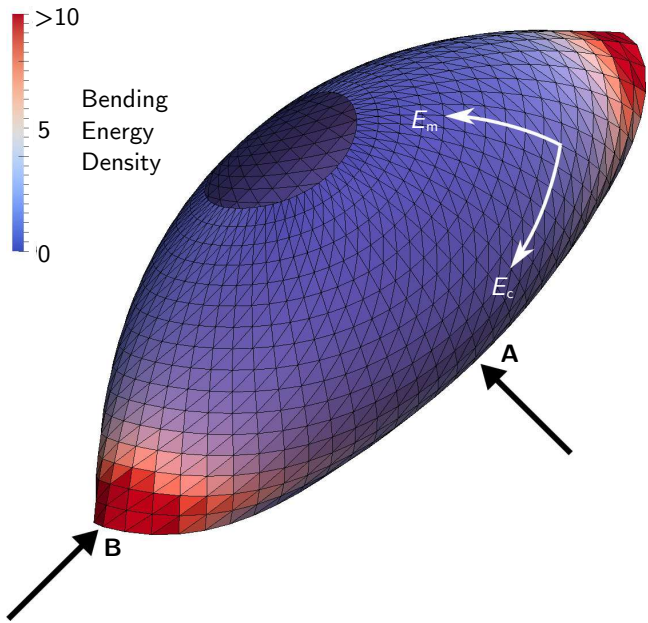}
\caption{Control mesh for deformed hemisphere ($\lambda = 0.5$) with its associated bending energy density after being loaded at \bs{A} and \bs{B} with pairwise opposite point loads}
\label{fig:hemi_result}
\end{figure}
This example is used to study how the shell performs when subjected to coupled stretching and bending stresses with large rigid body rotations. While this setup was originally developed for isotropic shells, it was modified to take into account the orthotropic case in Refs.\cite{Valdes2007,Valdes2009}. The problem consists of a hemisphere with a $18\degree$ hole at its north pole loaded by four equal point loads on its equator $90\degree$ from each other. These forces are diametrically opposite in direction, with a pair of tensile and compressive loads. The shell radius \is{R} is 10, with a thickness $h$ of 0.04, elastic modulus in the circumferential direction $E_\textrm{c}$ of $6.825 \times 10^7$, and a Poisson ratio $\nu_\textrm{c}$ of 0.3. The elastic modulus in the meridional direction $E_\textrm{m}$ and shear modulus $G_{cm}$ are given in Tab.(\ref{fig:load_cases}) for different material properties. For clarity, the degree of orthotropy is defined as $\lambda = E_\textrm{m}/E_\textrm{c}$.
\begin{table}[h]
\centering
\normalsize
\caption{Material properties for the pinched hemisphere benchmark \cite{Valdes2007,Valdes2009}}
\label{fig:load_cases}
\begin{tabular}{ c  c  c }
\toprule
$\lambda$ & $E_\mathrm{m}$ & $G_\mathrm{cm}$\\
\midrule
$1.0$ & $6.825 \times 10^7$ & $2.625 \times 10^7$ \\ 
$0.9$ & $6.143 \times 10^7$ & $2.518 \times 10^7$ \\
$0.5$ & $3.413 \times 10^7$ & $1.896 \times 10^7$ \\
$0.1$ & $6.825 \times 10^6$ & $5.884 \times 10^6$ \\
\bottomrule
\end{tabular}
\end{table}
The meshes were constructed by subdividing each triangular element into four until a fine enough discretization was obtained. Further data on the relevent isotropic case can be found in Refs.\cite{Simo1990,Sze2004}. Due to the small number of degrees of freedom (DOF) required for an accurate result and the lack of symmetry, the entire hemisphere was modelled. The hemispherical mesh consists of $16 \times 64$ subdivision elements as used by Ref.\cite{Vetter2013} with 3264 DOF (3 DOF per node). Fig.(\ref{fig:hemi_result}) illustrates two of the four diametrically opposite point loads (at \bs{A} and \bs{B}) on which an increasing load with a maximum of 100 is placed. Tab.(\ref{fig:pinch_results}) shows the maximum displacements of each material property scenario. The orthotropic subdivision shells (SD3R) are compared with quadratic shell elements with reduced integration (S8R) from the commercially available Abaqus software \cite{Abaqus.2011} and the rotation-free Basic Shell Triangle (BST) \cite{Valdes2007,Valdes2009} in Fig.(\ref{fig:load_comp}). Very good agreement with existing data leads us to conclude that the orthotropic subdivision shell is able to accurately model problems with large deformations and both stretching and bending with a substantially lower number of DOF than the BST and S8R elements found in literature.

\begin{figure*}[t] 
\centering
\includegraphics[width=12.9cm]{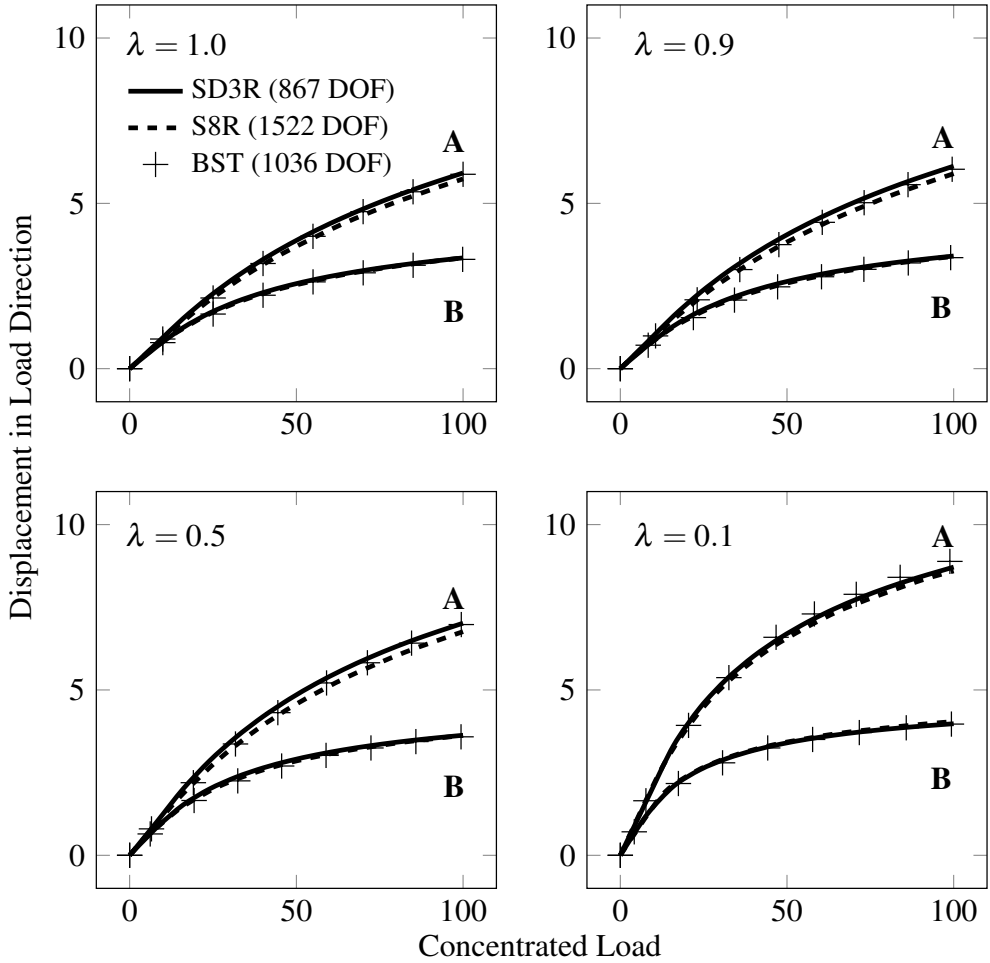}
\caption{Load-displacement curves for the deforming hemispherical shell with different degrees of orthotropy $\lambda$. The hemispherical SD3R mesh contains 3264 DOF for the entire hemisphere. For comparison with other methods, the number of DOF for a quarter of the hemisphere is $(16 \times 17) \times 3 = 867$ where 16 is the number of elements along the meridional axis. The S8R results were obtained using Abaqus \cite{Abaqus.2011} and the BST results from Refs.\cite{Valdes2007,Valdes2009}}\label{fig:load_comp}
\end{figure*}
\begin{table}[h]
\centering
\normalsize
\caption{Displacements of \bs{A} and \bs{B} at maximum load 100} \label{fig:pinch_results}
\begin{tabular}{ l  c  c  c  c }
\toprule
$\lambda$ & $1.0$ & $0.9$ & $0.5$ & $0.1$ \\
\midrule
A & 5.918 & 6.125 & 7.019 & 8.716 \\
B & 3.350 & 3.407 & 3.629 & 3.978 \\
\bottomrule
\end{tabular}
\end{table}
\subsection{Wrinkling of orthotropic sheets} \label{subsec:wrinkled}	
The wrinkling of sheets has been studied extensively with shell and membrane elements  \cite{Raible2005,Wong2006,Jarasjarungkiat2008,Flores2011}. The membrane element formulation needs to be augmented with a wrinkling model, but can be applied on coarser meshes, making it relatively computationally inexpensive. Subdivision shells have been used in the past to simulate the deformations of thin membranes \cite{Vetter2014}. To demonstrate the ability of orthotropic subdivision shell elements to reproduce the wrinkling phenomenon on coarse meshes, the reference shear test described in Ref.{\cite{Raible2005}} is simulated.
	
The test consists of a prestressed sheet sheared by displacement control. The objective is to find the critical shear displacement ($u_c$); a bifurcation point signaling the onset of wrinkling, and the maximum amplitude of the wrinkles at the maximum specified shear displacement. The 200 mm $\times$ 100 mm sheet is first prestressed with a displacement of 1 mm along its short axis, and then sheared by continuously displacing its upper edge by 10 mm along its long axis. The lower edge of the sheet remains pinned for the duration of the simulation (Fig.(\ref{fig:wrinkle})). The material is 0.2 mm thick and has a preferred orientation of $\alpha=30\degree$ to the \is{y}-axis. The material properties for the isotropic case are $E=600 \ \mathrm{N/mm^2}$ and $\nu=0.45$ \cite{Raible2005}, while for the orthotropic case $E_1 = E_2 = 106.6 \ \mathrm{N/mm^2}$, $\nu = 0.22$ and $G_{12} = 11.3 \ \mathrm{N/mm^2}$ \cite{Flores2011}. Meshes with resolution ranging from 8 $\times$ 4 (135 DOF) to 56 $\times$ 28 (4959 DOF) are simulated. 
	
\begin{figure}[h]
\centering
\includegraphics[width=8.4cm]{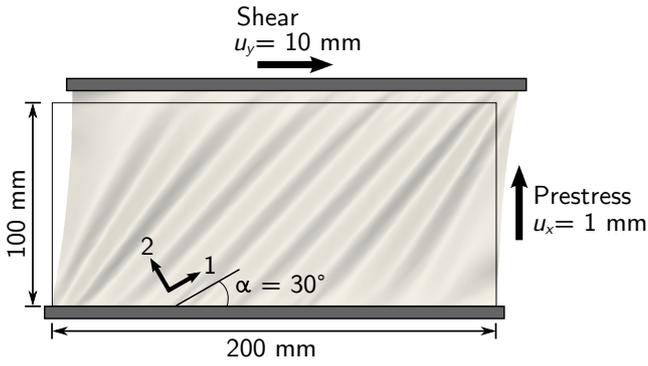}
\caption{Boundary and loading conditions of sheet to investigate wrinkling behavior \cite{Raible2005}}
\label{fig:wrinkle}
\end{figure}
In these simulations, a wrinkle is defined as a single fold in the sheet, which begins from the baseline $u_z = 0$, rises to its maximum value and back to $u_z = 0$. As seen from the results in Fig.({\ref{fig:convergence}}), the isotropic case shows a critical shear displacement of $u_c = 1.67$ mm, with a maximum wrinkle amplitude at a shear displacement of $u_y = 10$ mm of 1.9 mm which exactly corresponds to the value obtained with 3D membrane elements in Ref.\cite{Raible2005}. The number of wrinkles (14) for the maximum resolution of 56 $\times$ 28 elements (4959 DOF) also agrees with the existing data. In the orthotropic case, the critical shear displacement is $u_c = 1.81$ mm which compares favorably with Ref.\cite{Flores2011}. The maximum amplitude at a shear displacement of $u_y = 10$ mm is 2.2 mm, which is slightly larger than the 2 mm specified in the original study, while the sheet develops 10 wrinkles following the same profile. 
	
Strikingly, even a coarse mesh resolution of 16 $\times$ 8 elements (459 DOF) yields critical shear displacements of 1.68 mm and 1.81 mm for the isotropic and orthotropic case respectively, which are very small deviations from the high resolution meshes. The number of wrinkles converge more slowly for the isotropic case, requiring the 48 $\times$ 24 (3675 DOF) resolution to reach 14 wrinkles, while 10 wrinkles are already obtained for the 32 $\times$ 16 (1683 DOF) resolution in the orthotropic case. 

\begin{figure}[h]
\centering
\includegraphics[width=8.4cm]{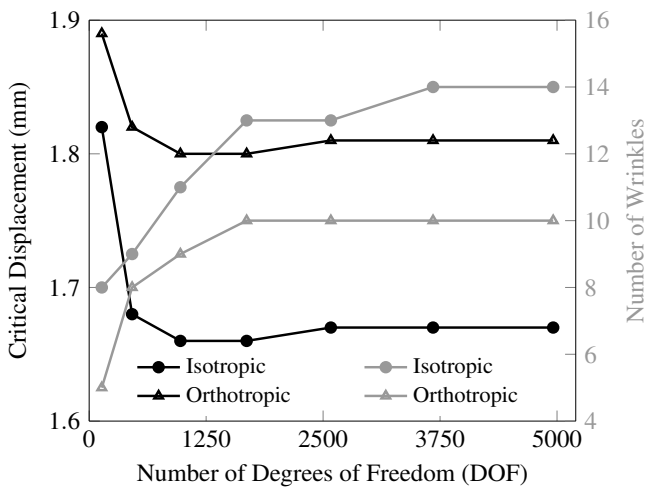}
\caption{Convergence of the critical shear displacement and number of wrinkles with mesh resolution}
\label{fig:convergence}
\end{figure}
\begin{figure}[h]
\centering
\includegraphics[width=8.4cm]{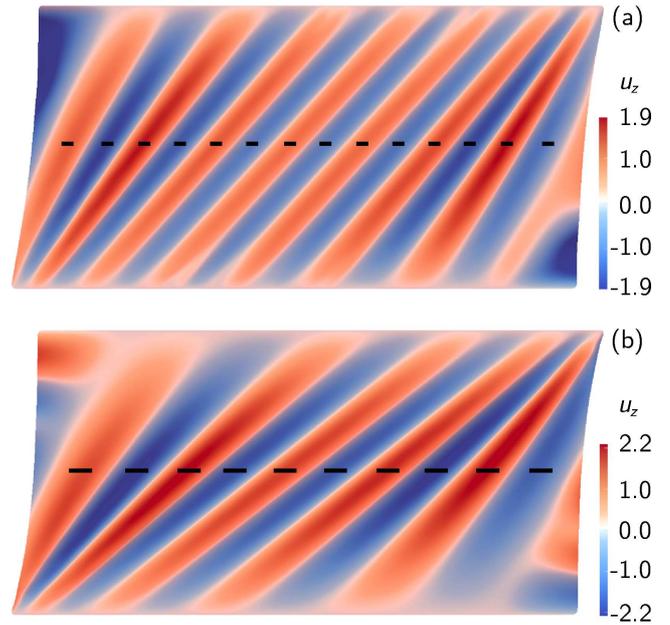}
\caption{Wrinkling of sheet with (a) isotropic \cite{Raible2005} and (b) orthotropic \cite{Flores2011} properties using a resolution of 56 $\times$ 28 elements with 4959 DOF. The displacement in the normal direction ($u_z$) is given in mm and the wrinkles are marked with black dashes}
\label{fig:contour}
\end{figure}
\section{Conclusions}
A method to simulate orthotropy in rotation-free shell finite elements was proposed. This approach transforms the derivatives of the shape functions to orthogonalize and align the basis for each triangular element. The implementation was performed using Kirchhoff--Love type subdivision shells. This approach requires negligible computational overhead as the transformation is only performed once in the undeformed configuration. The standard pinched hemispherical shell with $18\degree$ hole benchmark for orthotropic behaviour proved that the method is accurate. 

The combination of orthotropic material behaviour and the efficiency and robustness of subdivision finite elements has many advantages in simulating smooth membranes with arbitrary topologies. This implementation was then used to simulate the wrinkling of sheared orthotropic sheets where its efficiency was compared to existing methods. 

The shape function derivative transformation applied to traditional subdivision shell finite elements has been shown to be very well suited for the simulation of material anisotropies. The analysis of the sheared membrane showed that the wrinkling of an orthotropic sheet was accurate even for coarse meshes, with the added advantage that no specific wrinkling models are required to compensate for the lack of bending stiffness in membrane elements.   

\begin{acknowledgements}
The authors acknowledge support from the Research and Technology Development Project "MecanX: Physics-Based Models of Growing Plant Cells using Multi-Scale Sensor Feedback" granted under SystemsX.ch by the Swiss National Science Foundation, the Advanced Grant 319968-FlowCCS granted by the European Research Council (ERC) and from ETH Zurich by ETHIIRA Grant No. ETH-03 10-3.
\end{acknowledgements}

\bibliographystyle{spbasic}
\bibliography{Orthotropic}

\end{document}